\documentclass[10pt,twosided]{article}
\usepackage[tbtags]{amsmath}
\usepackage{epsfig,amstext,amssymb,amsthm,latexsym}
\allowdisplaybreaks[4]
\usepackage{color}

\usepackage{amssymb,color}
 
\newcommand{\EE}[1]{\mathbb{E}\left \{ #1 \right \}}

\usepackage{amssymb,color}

\definecolor{c20}{rgb}{0.,0.7,0.}
\definecolor{c30}{rgb}{0.,0.,1.}
\definecolor{c40}{rgb}{1,0.1,0.7}
\definecolor{c50}{rgb}{1,0,0}
\definecolor{c60}{rgb}{1,0.9,0.1}

\def\cE#1{\textcolor{c30}{#1}}
\def\cT#1{\textcolor{c50}{#1}}
\def\cE#1{#1}
\def\cT#1{#1}

\def\cH#1{#1}
\def\eT#1{#1}

\def\aH#1{#1}
\def\aT#1{#1}

\def\aE#1{\textcolor{c30}{#1}}
\def\aE#1{#1}

\def\peng#1{\textcolor{c30}{#1}}
\def\peng#1{#1}
\def\Tan#1{\textcolor{c40}{#1}}
\def\Tan#1{#1}

\def\enk#1{\textcolor{c20}{#1}}
\def\enk#1{#1}
\def\pzx#1{\textcolor{c30}{#1}}
\def\pzx#1{#1}

\def\tzq#1{\textcolor{c50}{#1}}
\def\tzq#1{#1}

\def\enke#1{\textcolor{c20}{#1}}
\def\enke#1{#1}

\usepackage{amssymb,color}

\def\kal#1{{\cal{ #1}}}

\newcommand{\ABs}[1]{ \biggl \lvert #1 \biggr \rvert}

\newcommand{\pk}[1]{\mathbb{P} \left \{ #1 \right\} }

\newcommand{\R}{\mathbb{R}}

\newcommand{\inr}{\in \R}

\newcommand{\inn}{\in \mathbb{N}}

\newcommand{\limit}[1]{\lim_{#1 \to   \infty}}

\newcommand{\BQN}{\begin{eqnarray}}
\newcommand{\EQN}{\end{eqnarray}}
\newcommand{\BQNY}{\begin{eqnarray*}}
\newcommand{\EQNY}{\end{eqnarray*}}

\newcommand{\BS}{\begin{sat}}
\newcommand{\ES}{\end{sat}}
\newcommand{\BT}{\begin{theo}}
\newcommand{\ET}{\end{theo}}
\newcommand{\BK}{\begin{korr}}
\newcommand{\EK}{\end{korr}}

\newcommand{\BD}{\begin{de}}
\newcommand{\ED}{\end{de}}

\newcommand{\BIT}{\begin{itemize}}
\newcommand{\EIT}{\end{itemize}}
\newcommand{\BDI}{\begin{description}}
\newcommand{\EDI}{\end{description}}

\newcommand{\BRM}{\begin{remarks}}
\newcommand{\ERM}{\end{remarks}}

\newcommand{\BEL}{\begin{lem}}
\newcommand{\EEL}{\end{lem}}

\newtheorem{theo}{Theorem}[section]

\newcommand{\COM}[1]{}

\newcommand{\QED}{\hfill $\Box$}

\def\Ha{{\cE{\cal H}_\alpha}}
\def\TH{\theta}
\def\tn{\aH{T(n)}}

\begin{document}
\baselineskip 15pt \setcounter{page}{1}
\title{\bf \Large  Asymptotics of Maxima of Strongly Dependent Gaussian Processes
}
\author{{\small  Zhongquan Tan\footnote{College of Mathematics Physics and Information Engineering, Jiaxing University, Jiaxing 314001, PR China},
\aH{Enkelejd Hashorva\footnote{Department of Actuarial Science, Faculty of Business and Economics, University of Lausanne, UNIL-Dorigny 1015 Lausanne, Switzerland}}, Zuoxiang Peng\footnote{School of Mathematics and Statistics, Southwest University, 400715 Chongqing,
China}}\\
}

\bigskip

\date{\today}
 \maketitle
 \baselineskip 15pt

{\bf Abstract:} \ Let $\{X_{n}(t), t\in[0,\infty)\}, n\in\mathbb{N}$
be \peng{a sequence of centered} dependent stationary Gaussian
process\aE{es}. \peng{The limit distribution of
$\sup_{t\in[0,T(n)]}|X_{n}(t)|$ is established as $r_{n}(t)$, the
correlation function of $\{X_{n}(t), t\in[0,\infty)\},
n\in\mathbb{N}$, satisfies the local and long range \enke{strong}
dependence conditions, which extends the results obtained by
Seleznjev (1991).}



\section{Introduction}

Let $\{X(t), t\in[0,\infty)\}$ be a standard \aE{(mean zero and unit
variance)} stationary Gaussian process with continuous sample
\peng{paths}, and let $\{r(t), t\ge 0\}$
denote its correlation function.  
Assume that the correlation function $r(t)$ \aE{of the process} satisfies
\begin{eqnarray}
\label{eq1.1}
r(t)=1-|t|^{\alpha}+o(|t|^{\alpha})\ \ \mbox{as}\ \ t\rightarrow 0,\quad  \mbox{and}\ \ r(t)<1\ \ \mbox{for}\ \ t>0
\end{eqnarray}
for some $\alpha\in(0, 2]$, and \aE{further} \peng{assume} 
\begin{eqnarray}
\label{eq1.2}
r(t)\log t\rightarrow 0,\ \ \mbox{as}\ \  t\rightarrow \infty.
\end{eqnarray}
\cE{For  the study of the asymptotic properties of the supremum of Gaussian processes the local condition \eqref{eq1.1} is a standard one, whereas the condition \eqref{eq1.2}
is the weak dependence condition, or the \enke{so-called} Berman's condition, see e.g., Piterbarg (1996). Under these two conditions on the correlation function $r(t)$}, it is well-known (see e.g., Leadbetter et al. (1983) \cE{or Berman (1992)}) that
\begin{eqnarray}
\label{eq1.3}
\lim_{T\to\infty}\enke{\sup_{x\inr}} \ABs{ \pk{a_{T}\left( \sup_{t\in[0,T]}X(t)-b_{T}\right)\leq x}
-\exp(-e^{-x})}=0,
\end{eqnarray}
where
\BQN \label{eq:AT}
a_{T}=\sqrt{2\log T}, \ \ b_{T}=\sqrt{2\log T}+\frac{\log( \Ha (2\pi)^{-1/2}  (2\log T)^{-1/2+1/\alpha})}{\sqrt{2\log
T}}.
\EQN
\cE{Here} $ \Ha $ denotes the Pickands constant defined by
$\aT{ \Ha =\lim_{\lambda\rightarrow\infty} \lambda^{-1} \Ha (\lambda)$,}
where
$$ \Ha (\lambda)=\EE{ \exp\left(\max_{t\in[0,\lambda]}\sqrt{2}B_{\alpha/2}(t)-t^{\alpha}\right)}$$
and $B_{\aE{\alpha}}$ is a fractional Brownian motion (a mean zero \aE{Gaussian process} with stationary increments such that
$\EE{B_{\aE{\alpha}}^{2}(t)}=|t|^{2\aE{\alpha}},\cE{t\inr} $). It is also well-known that
$0< \Ha <\infty$, see e.g., Berman (1992),  and Piterbarg (1996).\\
\cE{In this paper, the} following Pickands exact asymptotics plays a curial role in
deriving the limit relation of (\ref{eq1.3}). \aE{Specifically, for} some fixed constant $h>0$
\begin{eqnarray}
\label{eq1.4}
\pk{\sup_{t\in[0,h]}X(t)>u}&=&h\mu(u)(1+o(1)), \ \ \mbox{as}\ \ u\rightarrow\infty,
\end{eqnarray}
provided that \aE{the correlation function} $r(t)$ satisfies (\ref{eq1.1}) \aE{and}
\begin{eqnarray}
\label{eq1.5}
\mu(u)= \Ha u^{2/\alpha}\Psi(u),
\end{eqnarray}
where $\Psi(\cdot)$ is \cE{the survival function} of a standard
Gaussian random variable. \peng{For more details see  Leadbetter et
al.\ (1983) and Piterbarg (1996)}.
\enk{A correct proof of Pickand's theorem (see Pickands (1969)) was given in Piterbarg (1972); for the main properties of}
Pickands and related constants, see Adler (1990), Berman (1992), Shao (1996), Dieker (2005), D\c{e}bicki and Kisowski (2009) and
Albin and Choi (2010).

A uniform version of \cE{\eqref{eq1.4}} for stationary Gaussian
\cE{processes} has been \peng{established by} Seleznjev (1991),
where the author investigated the limit distribution of the error of
approximation of Gaussian stationary periodic \enke{processes} by random
trigonometric polynomials in the uniform metric. \cH{Next, we
formulate the aforementioned result.}

\bigskip

\textbf{Theorem A}. {\sl Let $\{X_{n}(t), t\in[0,\aH{\infty)}\}, n\in\mathbb{N}$ be standard stationary Gaussian processes with
a.s. continuous \aE{sample} paths and correlation function $r_{n}(t)$.  \aH{Let} $\tn>0,u_{n},n\ge 1$ be constants such that $\limit{n} \min(\tn,u_{n})=\infty$.} \enk{Suppose further that} 
\begin{itemize}
\item[(A1).] $r_{n}(t) = 1 -c_{n}|t|^{\alpha} + \varepsilon_{n}(t)|t|^{\alpha}$, $0<\alpha\leq 2$,
where $c_{n}\rightarrow1$ as $n\rightarrow\infty$ and
$\varepsilon_{n}(t)\rightarrow0$ as $t\rightarrow0$, uniformly in $n$\enk{;}
\item[({A2}).] for any $\varepsilon>0$, there exists $\gamma>0$ such that
$\sup\{|r_{n}(t)|, T\geq |t|\geq\varepsilon, n\in\mathbb{N}\}<\gamma<1$\enk{;}
\item[({A3}).] $r_{n}(t)\log(t)\rightarrow 0$ as $t\rightarrow\infty$, uniformly in $n$\enk{.}
\end{itemize}
\begin{itemize}
\item[(i).] If (A1) and (A2) hold, then for any fixed $h>0$ \enk{and $\mu(\cdot)$ defined in \eqref{eq1.5}}
\begin{eqnarray*}
\lim_{ n\to \infty} \frac{\pk{\sup_{t\in[0,h]}|X_{n}(t)|>u_{n}}}{2h\mu( u_{n})}&=&1.
\end{eqnarray*}
\item[(ii).] If additionally $\lim_{ n\to \infty} \tn\mu(u_{n})=\TH \in(0,\infty]$ and (A3) \enk{hold}, then
\begin{eqnarray*}
\lim_{ n\to \infty}\pk{\sup_{t\in[0,\tn]}|X_{n}(t)|\leq u_{n}} &=& e^{-2\TH},
\end{eqnarray*}
\enk{where we set $e^{-2\TH}=0$ if $\TH=\infty$.}
\item[(iii).]\cT{ If instead of Assumptions (A1)-(A3), the correlation functions $r_{n}(t)$ \aE{are such that} 
$$ 1-r_{n}(t) \leq |t|^{\alpha}, \quad  t\in[0, \tn],$$
\aE{with $\alpha\in (0,2]$ and  $\tn\geq T_{0}>0$} \enk{for all large $n$},  then
\begin{eqnarray*}
\lim_{ n\to \infty}\pk{\sup_{t\in[0, \tn]}|X_{n}(t)|\leq u_{n}} &=& 1,
\end{eqnarray*}
provided \peng{that} \aE{$\lim_{ n\to \infty} \tn\mu( u_{n})=0$.}}
\item[(iv).] \cT{Let $a_{\tn}, b_{\tn}$ \cE{be defined as in
\eqref{eq:AT}.} If (A1), (A2) and (A3) hold, then 
\begin{eqnarray*}
\enk{\lim_{ n\to \infty}\sup_{x\inr}} \ABs{\pk{ a_{\tn}\left(\sup_{t\in[0,\tn]}|X_{n}(t)|-b_{\tn}\right)\leq x}-\exp(-2e^{-x})}=0.
\end{eqnarray*}
}
\end{itemize}

\peng{The above result has been extended by Seleznjev (1996) to a
\cE{certain class} of non-stationary Gaussian process\cE{es}. For
further extensions and related studies, we refer to  H\"{u}sler
(1999),  H\"{u}sler et al. (2003) and Seleznjev (2006).}

With impetus from Seleznjev (1991), in this paper we present the
corresponding version of \peng{Theorem A for a sequence of} strongly
dependent stationary Gaussian process\cE{es} (see definition below).

The paper is organized as follows. Section 2 displays the main result, \cE{followed then by  Section  3 where} we present the proofs.

\section{Main Results}
In this section, \peng{we extend Theorem A  to a sequence of}
strongly dependent stationary Gaussian processes.  \cE{A sequence}
of \cE{ standard} stationary Gaussian process $\{X_{n}(t),
t\in[0,\cE{\infty)}\}, n\in\mathbb{N}$ is called strongly dependent
if the correlation function $r_{n}(t)$ satisfies \cE{one of the
following \enk{assumptions}}:
\begin{itemize}
\item[({B1}).] $r_{n}(t)\log t\rightarrow r\in (0, \infty)$ as  $t\rightarrow \infty$, uniformly in $n$;
\item[({B2}).] $r_{n}(t)\log t\rightarrow \infty$ as  $t\rightarrow \infty$, uniformly in $n$.
\end{itemize}
\cE{Indeed,} \Tan{Assumptions} \enk{(B1) and (B2)} are natural extensions of
\Tan{Assumption} (A3). For related studies on extremes for strongly
dependent Gaussian process, we refer to Mital and Ylvisaker (1975),
Piterbarg (1996), Ho and McCormick (1999) and Stamatovic and Stamatovic (2010).

\cE{Let \aE{in the following} $\varphi$ and \aT{$\Phi$ denote the
probability density  function and the \aE{distribution}  function}
of a standard Gaussian random variable \aE{$\kal{W}$}, respectively,
and set \BQN
\Lambda_{r}(x)
&=& \EE{ [\Lambda(x+r)]^{ e^{\sqrt{2r}\kal{W}}+e^{-\sqrt{2r}\kal{W}}}}, \quad x\inr,
\EQN
}
with $\Lambda(x)=\exp(-\exp(-x)), x\inr$ the unit Gumbel distribution function.

Next, we state our main results.

\BT \label{eq:thm1}
Let $\{X_{n}(t), t\in[0,\cE{\infty})\}, n\in\mathbb{N}$ be a standard stationary Gaussian processes with
a.s.\  continuous \aE{sample} paths and correlation function $r_{n}(t)$ satisfying (A1),(A2) and (B1).

(i). If $\aH{\lim_{n\to \infty}} T(n)\mu(u_{n})=\TH \in(0,\infty]$, then
\BQN\label{thm1:1}
\limit{n}\pk{\sup_{t\in[0,T(n)]}|X_{n}(t)|\leq u_{n}}
&=&\cE{\Lambda_{r}( - \log \theta)},
\EQN
\enk{where $\Lambda_{r}( - \log \theta)=:0$ if $\TH=\infty$}.\\
%
(ii). Let $a_{\tn}, b_{\tn}$ \cE{be defined as in \eqref{eq:AT}},
\peng{for $\Tan{x\inr}$ we have} \BQN \label{thm1:3}
\enk{\limit{n}\sup_{x\inr}}\ABs{ \pk{
a_{T(n)}\left(\sup_{t\in[0,T(n)]}|X_{n}(t)|-b_{T(n)}\right)\leq x}
- \cE{\Lambda_{r}( x)}}=0.
\EQN
\ET

\textbf{Remarks 2.1}. {\sl
(a) From the proof of Theorem 2.1, \enke{it follows that both \eqref{thm1:1} and \eqref{thm1:3} can be shown to hold also} for $r=0$, \cH{retrieving thus the result} of Theorem A.

(b) \cT{Assertion (iii) of Theorem A still holds under the conditions of Theorem 2.1.}
%
}

\BT \label{eq:thm2}
Let $\{X_{n}(t), t\in[0,\cE{\infty})\}, n\in\mathbb{N}$ be a standard stationary Gaussian processes with a.s.\ continuous sample paths
 and correlation function $r_{n}(t)$ satisfying (A1) \peng{with $0<\alpha\le 1$,} (A2) and (B2). Assume that $r_{n}(t)$ is convex for $t\geq 0$ and $r_n(t)=o(1)$ uniformly in $n$.
 If further $r_{n}(t)\log t$ is monotone for large $t$, then
with  $b_{\tn}$ \cE{as in \eqref{eq:AT}}, we have 
\BQN \label{thm1:3}
\enk{\limit{n} \sup_{x\in (0, \infty)}}\ABs{ \pk{ r_{n}^{-1/2}(\tn)\left(\sup_{t\in[0,T(n)]}|X_{n}(t)|-(1-r_{n}(\tn))^{1/2}b_{T(n)}\right)\leq x}
- 2\Phi(x)+1}=0.
\EQN
\ET

\textbf{Remarks 2.2}. {\sl Theorem 2.2 is a uniform version of Theorem 3.1 of Mittal and Ylvisaker (1975).}

\section{Further Results and Proofs}

\aE{We begin with some auxiliary lemmas needed for the proofs of Theorem 2.1 and 2.2. }

\cT{For given $\varepsilon> 0$, we divide interval $[0, \tn]$ onto intervals of length 1, and split each of them onto subintervals $I_{j}^{\varepsilon}$, $I_{j}$ of length $\varepsilon$, $1-\varepsilon$, $\tzq{j=1,2,\cdots,[\tn]},$ respectively,} where \tzq{$[x]$ denotes the integral part of $x$. It can be easily seen that a possible
remaining interval with length smaller than 1 plays no role in our consideration. We denote this interval with $J$.}

Let $\{X_{n}^{(i)}(t), t\geq 0\}$, $i=1,2,\cdots$ be independent copies of $\{X_{n}(t), \cE{t\ge 0}\}$ and
$\{\eta_{n}(t), t\geq 0\}$ be such that $\eta_{n}(t)=X_{n}^{(j)}(t)$ for \cT{$t\in I_{j}$}.
Let $\rho(\tn)\aH{:}=r/\log \tn$ and define
$$\xi_{n}(t)=(1-\rho(\tn))^{1/2}\eta_{n}(t)+\rho^{1/2}(\tn)\kal{W}, \ \ t\in \cup_{j=1}^{\tzq{[\tn]}} I_{j},$$
where $\kal{W}$ is a standard \aE{Gaussian random} variable
independent of \peng{$\{\eta_{n}(t),t\ge 0\}$}. Note that
$\{\xi_{n}(t), t\in \cup_{j=1}^{\tzq{[\tn]}} I_{j}\}$ is a standard
non-stationary Gaussian process \cE{with correlation function
$\varrho_{n}\cH{(\cdot,\cdot)}$ which is given by} \cT{\[
  \varrho_{n}(t,s)=\left\{
 \begin{array}{cc}
  {r_{n}(t,s)+(1-r_{n}(t,s))\rho(\tn)},    & t\in I_{j}, s\in I_{i}, i=j,\\
  {\rho(\tn)},    & t\in I_{j}, s\in I_{i}, i\neq j.
 \end{array}
  \right.
\]}

\peng{In the sequel, assume that $a, u_n,\aH{v_n}$ are positive
constants, and set}
$$q\enk{:=}q(u_{n})=au_{n}^{-2/\alpha}, \quad  \mu(u_{n}):= \Ha
u_{n}^{2/\alpha}\Psi(u_{n}), \quad \aH{\delta(a):=1-\frac{\Ha(a)}{\Ha}}.$$
Further, $C_1-C_6$ shall denote positive constants whose values may vary from place to place.

\textbf{Lemma 3.1}. {\sl \enk{If the Assumptions (A1) and (A2) hold, then for each interval $I$ of fixed length $h>0$}
\begin{eqnarray}
\label{eq3.2.1}
0\leq \pk{ \max_{jq\in I}|X_{n}(jq)|\leq u_{n}}-\pk{\sup_{s\in I}|X_{n}(s)|\leq u_{n}}\leq 2h\delta(a)\mu(u_{n})+o(\mu(u_{n}))
\end{eqnarray}
and
\begin{eqnarray}
\label{eq3.2.2}
0\leq \pk{\max_{jq\in I}X_{n}(jq)\leq u_{n}}-\pk{\sup_{s\in I}X_{n}(s)\leq u_{n}}\leq h\delta(a)\mu(u_{n})+o(\mu(u_{n})),
\end{eqnarray}
where $\delta(a)\rightarrow 0$ as $a\downarrow 0$.
}

\begin{proof} \cE{Both claims above are established in} the proof of Theorem 1 of Seleznjev (1991).
\end{proof}

\textbf{Lemma 3.2}. {\sl  Suppose that (A1) and (A2) hold. If  $T(n)\mu(u_{n})=O(1)$ and $T(n)\mu(v_{n})=O(1)$, then
\begin{eqnarray}
\label{eq3.3.1}
\pk{\sup_{s\in[0,\tn]}|X_{n}(s)|\leq u_{n}}-\pk{\sup_{s\in\cup I_{j}}|X_{n}(s)|\leq u_{n}} &\to& 0
\end{eqnarray}
and
\begin{eqnarray}
\label{eq3.3.2}
\pk{-v_{n}\leq \inf_{s\in[0,1]}X_{n}(s), \sup_{s\in[0,1]}X_{n}(s)\leq u_{n}}
 -\pk{-v_{n}\leq \inf_{s\in I_{1}}X_{n}(s), \sup_{s\in I_{1}}X_{n}(s)\leq u_{n}} &\to& 0
\end{eqnarray}
as $n\rightarrow\infty$ and $\varepsilon\downarrow 0$.
}

\begin{proof} \cE{By the} stationarity of $\{X_{n}(t), t\in[0, \tn]\}$ and Theorem A \cH{(i)}  \cE{we obtain }
\begin{eqnarray*}
&&\left|\pk{\sup_{s\in[0,\tn]}|X_{n}(s)|\leq u_{n}}-\pk{\sup_{s\in\cup I_{j}}|X_{n}(s)|\leq u_{n}}\right|\\
&&\leq \sum_{j=1}^{\tzq{[\tn]}}\pk{\max_{s\in I_{j}^{\varepsilon}}|X_{n}(s)|>u_{n}}+\tzq{\pk{\max_{s\in J}|X_{n}(s)|>u_{n}}}\\
&&\leq2(\tzq{[\tn]\varepsilon+1}) \mu(u_{n})(1+o(1))\\
&&=O(1)\varepsilon(1+o(1))\\
&&\rightarrow 0
\end{eqnarray*}
as $u\rightarrow\infty$ and $\varepsilon\downarrow 0$, which completes the proof of (\ref{eq3.3.1}).
Note \cE{in passing} that
\BQNY
\lefteqn{\left| \pk{-v_{n}\leq \inf_{s\in[0,1]}X_{n}(s), \sup_{s\in[0,1]}X_{n}(s)\leq u_{n}}
 -\pk{-v_{n}\leq \inf_{s\in I_{1}}X_{n}(s), \sup_{s\in I_{1}}X_{n}(s)\leq u_{n}}\right|}\\
&\leq &\left|\pk{\sup_{s\in[0,1]}X_{n}(s)\leq u_{n}}-\pk{\sup_{s\in I_{1}}X_{n}(s)\leq u_{n}}\right|\\
&+& \left|\pk{\inf_{s\in[0,1]}X_{n}(s)\geq -v_{n}}-\pk{\inf_{s\in I_{1}}X_{n}(s)\geq -v_{n}}\right|.
\EQNY
The proof of (\ref{eq3.3.2})  is similar to that of (\ref{eq3.3.1}), and therefore omitted.
\end{proof}

\textbf{Lemma 3.3}. {\sl  \enk{Under the assumptions of Lemma 3.2 we have} 
\begin{eqnarray}
\label{eq3.4.1}
\pk{\sup_{s\in\cup I_{j}}|X_{n}(s)|\leq u_{n}}-\pk{\max_{kq\in\cup I_{j}}|X_{n}(kq)|\leq u_{n}} &\rightarrow & 0
\end{eqnarray}
and
\begin{eqnarray}
\label{eq3.4.2}
\pk{-v_{n}\leq \inf_{s\in I_{1}}X_{n}(s), \sup_{s\in I_{1}}X_{n}(s)\leq u_{n}}-
\pk{-v_{n}\leq \min_{kq\in I_{1}}X_{n}(kq), \max_{kq\in I_{1}}X_{n}(kq)\leq u_{n}} &\rightarrow &0
\end{eqnarray}
as $n\rightarrow\infty$ and $a\downarrow 0$.
}

\begin{proof} \cE{By Lemma 3.2}
\begin{eqnarray*}
&&\left|\pk{\sup_{s\in\cup I_{j}}|X_{n}(s)|\leq u_{n}}-\pk{\sup_{kq\in\cup I_{j}}|X_{n}(kq)|\leq u_{n}}\right|\\
&&\leq \peng{T(n)}\max_{j}\left(\pk{\max_{kq\in I_{j}}|X_{n}(kq)|\leq u_{n}}-\pk{\sup_{s\in I_{j}}|X_{n}(s)|\leq u_{n}}\right)\\
&&\leq 2(1-\varepsilon) \tzq{[\tn]} \mu(u_{n})\delta(a)+\tn o(\mu(u_{n}))\\
&&= 2(1-\varepsilon)O(1)\delta(a)+o(1)\\
&&\rightarrow 0
\end{eqnarray*}
as $n\rightarrow\infty$ and $a\downarrow 0$. Hence the first claim follows. 
Note that
\begin{eqnarray*}
\lefteqn{\left|\pk{-v_{n}\leq \inf_{s\in I_{1}}X_{n}(s), \sup_{s\in I_{1}}X_{n}(s)\leq u_{n}}-\pk{-v_{n}\leq \min_{kq\in I_{1}}X_{n}(kq), \max_{kq\in I_{1}}X_{n}(kq)\leq u_{n}}\right|}\\
&\leq & \left|\pk{\max_{kq\in I_{1}}X_{n}(kq)\leq u_{n}}-\pk{\sup_{s\in I_{1}}X_{n}(s)\leq u_{n}}\right|
+\left|\pk{\min_{kq\in I_{1}}X_{n}(kq)\geq -v_{n}}-\pk{\inf_{s\in I_{1}}X_{n}(s)\geq -v_{n}}\right|.
\end{eqnarray*}
\cE{We omit the proof of} (\ref{eq3.4.2})  since it is similar to that of (\ref{eq3.4.1}).
\end{proof}

\textbf{Lemma 3.4}. {\sl  Suppose that (A1),(A2) and (B1) hold. If  $T(n)\mu(u_{n})=O(1)$, then
\begin{eqnarray}
\label{eq3.1}
\lim_{n  \to \infty}\left|\pk{\max_{kq\in\cup I_{j}}|X_{n}(kq)|\leq u_{n}}-\pk{\max_{kq\in \cup I_{j}}|\xi_{n}(kq)|\leq u_{n}}\right|&\aH{=}& 0.
\end{eqnarray}
}

\def\rh{\aH{r^{(h)}}}
\begin{proof}
Applying \Tan{ the generalized Berman inequality
(cf. Theorem 1.2 of Piterbarg (1996))}, we have \aH{(set next
$T:=\tn$)}
\begin{eqnarray}
\label{eq3.1.1}
&&\left|\pk{\max_{kq\in\cup I_{j}}|X_{n}(kq)|\leq u_{n}}-\pk{\max_{kq\in \cup I_{j}}|\xi_{n}(kq)|\leq u_{n}}\right|\nonumber\\
&&\leq \sum_{ \pzx{kq\in I_{i},\; lq\in I_{j}}}
\pzx{\frac{4}{2\pi}}|r_{n}(kq,lq)-\varrho_{n}(kq,lq)|\int_{0}^{1}\frac{1}{\sqrt{1-
\rh(kq,lq)}}\exp\left(-\frac{u_{n}^{2}}{1+\rh
(kq,lq)}\right)dh\nonumber\\
&&\pzx{\le\sum_{kq\in I_{i},\; lq\in I_{i},\atop i\in\{1,2,\cdots,[T(n)]\}}
\mathbb{A}(n,k,l,q)+ \sum_{kq\in I_{i},\; lq\in I_{j},i\neq j\atop i,j\in\{1,2,\cdots,[T(n)]\}}
\mathbb{A}(n,k,l,q)},
\end{eqnarray}
where $\varphi(x,y,\rh)$ is a Gaussian two-dimensional density with the covariance $\rh$, the variance equal to one
and zero mean and
$$\rh(kq,lq)=hr_{n}(kq,lq)+(1-h)\varrho_{n}(kq,lq), \quad \aH{h \in [0,1]}.$$
In the following part of the proof, let
$\varpi_{n}(kq,lq)=\max\{|r_{n}(kq,lq)|,|\varrho_{n}(kq,lq)|\}$ and
$\vartheta(t)=\sup_{t<|kq-lq|\leq T}\{\varpi_{n}(kq,lq)\}$. By
Assumption (A2) and \enke{the definition of} $\varrho_{n}(t,s)$, we
have \peng{$\vartheta(\varepsilon)=\sup_{\varepsilon<|kq-lq|\leq
T}\{\varpi_{n}(kq,lq); n\in\mathbb{N}\}<1$ for sufficiently large
\aH{$T$}}. Further, let $\beta$ be such that \peng{
$0<\beta<\frac{1-\vartheta(\varepsilon)}{1+\vartheta(\varepsilon)}$ for all sufficiently large $T$}.\\
\pzx{\enke{Next}, \cH{we} estimate the upper bound of}
(\ref{eq3.1.1}) in the case that $kq$ and $lq$ belong to the same
interval $I$. Note that in this case,
$\varrho_{n}(kq,lq)=r_{n}(kq,lq)+(1-r_{n}(kq,lq))\rho(T)\sim
r_{n}(kq,lq)$ for sufficiently large $T$. Split \pzx{the first term
of} (\ref{eq3.1.1}) into two parts as\pzx{
\begin{eqnarray}
\label{eq3.1.A} \sum_{kq\in I_{i},\; lq\in I_{i},i\in\{1,2,\cdots,[T(n)]\}\atop 0<|kq-lq|\leq \varepsilon
}\mathbb{A}(n,k,l,q) + \sum_{kq\in I_{i},\; lq\in I_{i},i\in\{1,2,\cdots,[T(n)]\}\atop
\varepsilon<|kq-lq|\leq 1-\varepsilon }
\mathbb{A}(n,k,l,q)=:J_{n1}+J_{n2}.
\end{eqnarray}}
The Assumption (A1) implies for all
$|t|\leq\varepsilon<2^{-1/\alpha}$
$$1-r_{n}(t)\leq 2|t|^{\alpha}.$$
From the assumption that $\aH{T\mu(u_{n})}= T(n)\mu(u_{n})=O(1)$,
\peng{we have
\begin{eqnarray}
\label{eq3.1.B}
u_{n}\sim (2\log T)^{1/2},\quad
e^{-\frac{u_{n}^{2}}{2}}\sim
(2\pi)^{1/2}H_{\alpha}^{-1}u_{n}^{1-2/\alpha}T^{-1}O(1).
\end{eqnarray}
Consequently}, \enk{with} $q:=au_{n}^{-2/\alpha}\sim a(\log
T)^{-1/\alpha}$ we obtain
\begin{eqnarray}
\label{eq3.1.C}
J_{n1}&\leq&C_{1}\sum_{kq\in I_{i},\; lq\in I_{i},i\in\{1,2,\cdots,[T(n)]\}\atop 0<|kq-lq|\leq \varepsilon}|r_{n}(kq,lq)-\peng{\varrho_{n}}(kq,lq)|\frac{1}{\sqrt{1-\varrho_{n}(kq,lq)}}\exp\left(-\frac{u_{n}^{2}}{1+\varrho_{n}(kq,lq)}\right)\nonumber\\
&\leq & C_{1}\sum_{kq\in I_{i},\; lq\in I_{i},i\in\{1,2,\cdots,[T(n)]\} \atop 0<|kq-lq|\leq \varepsilon}|(1-r_{n}(kq,lq))\rho(T)|\frac{1}{\sqrt{1-r_{n}(kq,lq)}}\exp\left(-\frac{u_{n}^{2}}{1+r_{n}(kq,lq)}\right)\nonumber\\
&\leq & C_{1}\rho(T)\sum_{kq\in I_{i},\; lq\in I_{i},i\in\{1,2,\cdots,[T(n)]\} \atop 0<|kq-lq|\leq \varepsilon}\sqrt{1-r_{n}(kq,lq)}\exp\left(-\frac{u_{n}^{2}}{1+r_{n}(kq,lq)}\right)\nonumber\\
&\leq & C_{1}\rho(T)\frac{T}{q}\sum_{0<kq\leq \varepsilon}\sqrt{1-r_{n}(kq)}\exp\left(-\frac{u_{n}^{2}}{2}\right)\exp\left(-\frac{(1-r_{n}(kq))u_{n}^{2}}{2(1+r_{n}(kq))}\right)\nonumber\\
&\leq & C_{1} \rho(T)\frac{T}{q}T^{-1}(\log T)^{1/2-1/\alpha}\sum_{0<kq\leq \varepsilon}(kq)^{\alpha/2}\exp\left(-\frac{1}{8}|kq|^{\alpha}\right)\nonumber\\
&\leq & C_{1} (\log T)^{-1/2},
\end{eqnarray}
which implies $\cH{\limit{n} J_{n1}=0}.$ By \peng{\eqref{eq3.1.B} for large $T$ we have
\begin{eqnarray}
\label{eq3.1.D}
J_{n2}&\leq & C_{2} \sum_{kq\in I_{i},\; lq\in I_{i},i\in\{1,2,\cdots,[T(n)]\}\atop \varepsilon<|kq-lq|\leq 1-\varepsilon }|r_{n}(kq,lq)-\peng{\varrho_{n}}(kq,lq)|\exp\left(-\frac{u_{n}^{2}}{1+\varpi_{n}(kq,lq)}\right)\nonumber\\
&\leq & C_{2} \sum_{kq\in I_{i},\; lq\in I_{i},i\in\{1,2,\cdots,[T(n)]\}\atop \varepsilon<|kq-lq|\leq 1-\varepsilon }\exp\left(-\frac{u_{n}^{2}}{1+\vartheta(\varepsilon)}\right)\nonumber\\
&\leq & C_{2}\frac{T}{q}\sum_{\varepsilon<kq\leq 1-\varepsilon}\exp\left(-\frac{u_{n}^{2}}{1+\vartheta(\varepsilon)}\right)\nonumber\\
&\leq & C_{2} \frac{T}{q^{2}}\left(\exp\left(-\frac{u_{n}^{2}}{2}\right)\right)^{\frac{2}{1+\vartheta(\varepsilon)}}\nonumber\\
&\leq & C_{2}T^{-\frac{1-\vartheta(\varepsilon)}{1+\vartheta(\varepsilon)}}(\log T)^{\frac{2\vartheta(\varepsilon)+\alpha}{\alpha(1+\vartheta(\varepsilon))}}.
\end{eqnarray}
Hence since $\vartheta(\varepsilon)<1$,} then
$\limit{n} J_{n2}=0$. \\
\enk{We continue with an estimate for}  \pzx{the upper bound of}
(\ref{eq3.1.1}) where $kq\in I_{i}$ and $lq\in I_{j}$, $i\neq j$.
Note that in this case, the distance between any two intervals
$I_{i}$ and $I_{j}$ is large than $\varepsilon$. Split the
\pzx{second term} of (\ref{eq3.1.1}) \enke{as} \pzx{
\begin{eqnarray}
\label{eq3.1.2} \sum_{kq\in I_{i},\; lq\in I_{j},i\neq j\in\{1,2,\cdots,[T(n)]\}\atop\varepsilon<|kq-lq|\leq
T^{\beta}}\mathbb{A}(n,k,l,q)+\sum_{kq\in I_{i},\; lq\in I_{j},i\neq j\in\{1,2,\cdots,[T(n)]\}\atop T^{\beta}<|kq-lq|\leq
T}\mathbb{A}(n,k,l,q)=:I_{n1}+I_{n2}.
\end{eqnarray}}
\enk{Similarly to the derivation} of (\ref{eq3.1.D}),  we  \peng{have
\begin{eqnarray}
\label{eq3.1.3}
I_{n1}&\leq&C_{3}\sum_{kq\in I_{i},\; lq\in I_{j},i\neq j\in\{1,2,\cdots,[T(n)]\}\atop\varepsilon<|kq-lq|\leq T^{\beta}}|r_{n}(kq,lq)-\peng{\varrho_{n}}(kq,lq)|\exp\left(-\frac{u_{n}^{2}}{1+\varpi_{n}(kq,lq)}\right)\nonumber\\
&\leq &C_{3}\sum_{kq\in I_{i},\; lq\in I_{j},i\neq j\in\{1,2,\cdots,[T(n)]\}\atop\varepsilon<|kq-lq|\leq T^{\beta}}\exp\left(-\frac{u_{n}^{2}}{1+\vartheta(\varepsilon)}\right)\nonumber\\
&\leq & C_{3}\frac{T}{q}\sum_{\varepsilon<kq\leq T^{\beta}}\exp\left(-\frac{u_{n}^{2}}{1+\vartheta(\varepsilon)}\right)\nonumber\\
&\leq & C_{3} \frac{T^{1+\beta}}{q^{2}}\left(\exp\left(-\frac{u_{n}^{2}}{2}\right)\right)^{\frac{2}{1+\vartheta(\varepsilon)}}\nonumber\\
&\leq & C_{3}T^{\beta-\frac{1-\vartheta(\varepsilon)}{1+\vartheta(\varepsilon)}}(\log T)^{\frac{2\vartheta(\varepsilon)+\alpha}{\alpha(1+\vartheta(\varepsilon))}}.
\end{eqnarray}}
Thus, $\aH{\lim_{n\to \infty}}I_{n1}=0$, since \peng{$\beta<\frac{1-\vartheta(\varepsilon)}{1+\vartheta(\varepsilon)}$}.
\enk{Further}, \cH{Assumption} (B1) implies \peng{that there \enk{exists a}
\enke{positive} constant $K$ such that $\varpi_{n}(kq)\leq K/\log T^{\beta}$ for
$kq>T^{\beta}$}. Using (\ref{eq3.1.B}) again, for
\peng{$q=au_{n}^{-2/\alpha}\sim a(\log T)^{-1/\alpha}$ we have
\begin{eqnarray*}
\frac{T^{2}}{q^{2}\log
T}\exp\left(-\frac{u_{n}^{2}}{1+\vartheta(T^{\beta})}\right)
&\leq &\frac{T^{2}}{q^{2}\log T}\exp\left(-\frac{u_{n}^{2}}{1+K/\log T^{\beta}}\right)\\
&\leq &C_{4}\exp\left(\frac{2K\log T}{K+\beta\log T}-\left(1-2/\alpha\right)\frac{K\log\log T}{K+\beta\log T}\right)\\
&=&O(1).
\end{eqnarray*}}
\cE{Hence,} following the argument of the proof of Lemma 6.4.1 of
Leadbetter et al. (1983) we may further write\peng{
\begin{eqnarray}
\label{eq3.1.4}
I_{n2}&\leq&C_{5}\sum_{kq\in I_{i},\; lq\in I_{j},i\neq j\in\{1,2,\cdots,[T(n)]\}\atop T^{\beta}<|kq-lq|\leq T}|r_{n}(kq,lq)-\varrho_{n}(kq,lq)|\exp\left(-\frac{u_{n}^{2}}{1+\varpi_{n}(kq,lq)}\right)\nonumber\\
&\leq&C_{5}\sum_{kq\in I_{i},\; lq\in I_{j},i\neq j\in\{1,2,\cdots,[T(n)]\}\atop T^{\beta}<|kq-lq|\leq T}|r_{n}(kq,lq)-\rho(T)|\exp\left(-\frac{u_{n}^{2}}{1+\vartheta(T^{\beta})}\right)\nonumber\\
&=&C_{5} \frac{q\log T}{T}\cE{\sum_{T^{\beta}<kq\leq T}|r_{n}(kq)-\rho(T)|}\frac{T^{2}}{q^{2}\log T}\exp\left(-\frac{u_{n}^{2}}{1+\vartheta(T^{\beta})}\right)\nonumber\\
&\leq& C_{5} \frac{q\log T}{T}\sum_{T^{\beta}<kq\leq T }|r_{n}(kq)-\rho(T)|\nonumber\\
&\leq&  \cE{C_5}\frac{q}{\beta T}\sum_{T^{\beta}<kq\leq T
}|r_{n}(kq)\log kq-r|+C_6r\frac{q}{T}\sum_{T^{\beta}<kq\leq T
}|1-\frac{\log T}{\log kq}|.
\end{eqnarray}}
\cE{By Assumption (B1),} the first term of the right hand-side of (\ref{eq3.1.4}) tends to 0.
Furthermore, the second term \cE{therein} also tends to 0, \cE{which follows}  by an integral estimate as in the proof of Lemma 6.4.1 of Leadbetter et al. (1983).
\cE{Consequently, the proof is established by \eT{(\ref{eq3.1.1})-(\ref{eq3.1.A}) and (\ref{eq3.1.C})-(\ref{eq3.1.4}).} }
\end{proof}

\textbf{Lemma 3.5}. {\sl  Suppose that (A1) and (A2) hold. If  $T(n)\mu(u_{n})=O(1)$ and $T(n)\mu(v_{n})=O(1)$, then
\begin{eqnarray}
\label{eq3.5.1}
\pk{\sup_{s\in[0,1]}X_{n}(s)> u_{n}, \inf_{s\in[0,1]}X_{n}(s)< -v_{n}}=o(\mu(u_{n})+\mu(v_{n})),  \quad \enk{n\rightarrow\infty.}
\end{eqnarray}
}

\begin{proof}
\enk{The proof is similar to that of} Lemma 11.1.4 in Leadbetter et al. (1983).
\end{proof}

\bigskip

{\it Proof of Theorem} 2.1. We only prove case (i), since case (ii) is a special case of (i).\\
(1). Case $\theta\in(0,\infty)$.
\enk{The definition} of $\{\xi_{n}(t), t\in \cup_{j=1}^{\tzq{[\tn]}} I_{j}\}$ \enk{implies}
\BQN
\label{eq3.2}
\pk{\max_{kq\in \cup I_{j}}|\xi_{n}(kq)|\leq u_{n}}&=& \pk{\max_{kq\in \cup I_{j}}|(1-\rho(\tn))^{1/2}\eta_{n}(kq)+\rho^{1/2}(\tn)\peng{\kal{W}}|\leq u_{n}}\nonumber\\
&=&\pk{-u_{n}\leq(1-\rho(\tn))^{1/2}\eta_{n}(kq)+\rho^{1/2}(\tn)\peng{\kal{W}}\leq u_{n},kq\in \cup I_{j}}\nonumber\\
&=&\int_{-\infty}^{+\infty}\pk{\frac{-u_{n}-\rho^{1/2}(\tn)z}{(1-\rho(\tn))^{1/2}}\leq\eta_{n}(kq)\leq
\frac{u_{n}-\rho^{1/2}(\tn)z}{(1-\rho(\tn))^{1/2}},kq\in \cup
I_{j}}\varphi(z)\, dz. \EQN
\enk{Since as $n\rightarrow\infty$} 
\begin{eqnarray*}
u_{n}^{(z)}&:=&\frac{u_{n}-\rho^{1/2}(\tn)z}{(1-\rho(\tn))^{1/2}}
=u_{n}+\frac{r-\sqrt{2r}z}{u_{n}}+o(u_{n}^{-1})
\end{eqnarray*}
and
\begin{eqnarray*}
v_{n}^{(z)}&:=&\frac{u_{n}+\rho^{1/2}(\tn)z}{(1-\rho(\tn))^{1/2}}
=u_{n}+\frac{r+\sqrt{2r}z}{u_{n}}+o(u_{n}^{-1}).
\end{eqnarray*}
\pzx{So, the assumption}
$\limit{n}T(n)\mu(u_{n})=\TH \in(0,\infty)$ \enk{implies} that
\begin{eqnarray}
\label{eq3.3}
\cE{\limit{n}}\tn \mu(u_{n}^{(z)})= \TH e^{-r+\sqrt{2r}z},\ \ \limit{n} \tn\mu(v_{n}^{(z)})= \TH e^{-r-\sqrt{2r}z}.
\end{eqnarray}
\cE{Next,} by the definition of $\{\eta_{n}(t),t\geq 0\}$, (\ref{eq3.3.2}), (\ref{eq3.4.2}) and (\ref{eq3.3}) we have
\begin{eqnarray}
\pk{-v_{n}^{(z)}\leq\eta_{n}(kq)\leq u_{n}^{(z)},kq\in \cup I_{j}}&=&\prod_{j=1}^{\tzq{[\tn]}}\pk{-v_{n}^{(z)}\leq X_{n}^{(j)}(kq)\leq u_{n}^{(z)}, kq\in I_{j}}\nonumber\\
&=&\pk{-v_{n}^{(z)}\leq X_{n}(kq)\leq u_{n}^{(z)}, kq\in I_{1})}^{\tzq{[\tn]}}\nonumber\\
&=&\pk{-v_{n}^{(z)}\leq X_{n}(t)\leq u_{n}^{(z)}, t\in I_{1}}^{\tzq{[\tn]}}(1+o(1))\nonumber\\
&=&\pk{-v_{n}^{(z)}\leq X_{n}(t)\leq u_{n}^{(z)},t\in [0,1]}^{\tzq{[\tn]}}(1+o(1))\nonumber\\
&=&\Biggl(
1-\pk{\inf_{s\in [0,1]}X_{n}(s)<-v_{n}^{(z)}}-\pk{\sup_{s\in [0,1]}X_{n}(t)> u_{n}^{(z)}} \notag\\
&& +\pk{\inf_{s\in [0,1]}X_{n}(s)<-v_{n}^{(z)}, \sup_{s\in [0,1]}X_{n}(t)> u_{n}^{(z)}}\Biggr)^{\tzq{[\tn]}}(1+o(1))
\end{eqnarray}
as $n\rightarrow\infty$.
\cE{In the light of} Theorem A(i) and Lemma 3.5
\BQNY
\pk{-v_{n}^{(z)}\leq\eta_{n}(kq)\leq u_{n}^{(z)},kq\in \cup I_{j}}&=&\left(1-\mu(u_{n}^{(z)})-\mu(v_{n}^{(z)})+o(\mu(u_{n}^{(z)})+\mu(v_{n}^{(z)}))\right)^{\tzq{[\tn]}}(1+o(1))\nonumber\\
&=&\left(1-\frac{\TH e^{-(r-\sqrt{2r}z)}+\TH e^{-(r+\sqrt{2r}z)}}{T(n)}+o\left(\frac{1}{\tn}\right)\right)^{\tzq{[\tn]}}(1+o(1))\nonumber\\
&=&\exp\left(-\TH e^{-(r-\sqrt{2r}z)}-\TH  e^{-(r+\sqrt{2r}z)}\right)(1+o(1))
\EQNY
as $n\rightarrow\infty$. Combining \cE{the last result} with (\ref{eq3.1}),(\ref{eq3.2})
\cE{and applying the} dominated convergence theorem we have
\BQNY
\limit{n} \pk{\max_{kq\in\cup I_{j}}|X_{n}(kq)|\leq u_{n}}
&=& \int_{-\infty}^{+\infty}\exp\left(-\TH e^{-(r-\sqrt{2r}z)}-\TH  e^{-(r+\sqrt{2r}z)}\right)\varphi(z)\, dz.
\EQNY
\cE{Consequently, the proof follows} \cE{by utilising further} (\ref{eq3.3.1}), (\ref{eq3.4.1}) and (\ref{eq3.1}).\\
\cT{(2). Case $\theta=\infty$.  From the definition of $\mu(\cdot)$, we know that for arbitrarily large $\theta'<\infty$, there exist a real sequence $v_{n}$ such that $\aH{\lim_{n\to \infty}} n\mu(v_{n})=\theta'$.
Clearly, for $n$ sufficient large, $u_{n}\leq v_{n}$, \aH{hence}
$$\pk{\sup_{t\in[0,T(n)]}|X_{n}(t)|\leq u_{n}}\leq\pk{\sup_{t\in[0,T(n)]}|X_{n}(t)|\leq v_{n}}\rightarrow \Lambda_{r}( - \log \theta'), \quad \aH{n\to \infty.}$$
Since this holds for arbitrarily large $\theta'<\infty$, by letting $\theta'\rightarrow \infty$ we see that
$$\lim_{n \to \infty}\pk{\sup_{t\in[0,T(n)]}|X_{n}(t)|\leq u_{n}}=0,$$
which completes the proof.}
\QED

\aE{For the proof of Theorem 2.2 we need a result which is
formulated in the next lemma. By Polya's criterion (see e.g., (3.10)
\aE{in} Durrett 2004) if we assume the convexity of the correlation
functions} $r_{n}(t)$ \pzx{(hence $0<\alpha\le 1$, cf. Theorem
 3.1 of Mittal and Ylvisaker (1975))}, then there exists a separable
standard stationary Gaussian process $Y_{n}(t),n\inn$ with correlation
function
$$\rho_{n,\tn}(t)= \frac{r_{n}(t)-r_{n}(\tn)}{1-r_{n}(\tn)}, \quad \ \mbox{for}\ \ \ t\leq \tn.$$
Let
$$M_{\tn}(Y)=\max_{0\leq t\leq \tn}Y_{n}(t),\quad   M_{\tn}(-Y)=\max_{0\leq t\leq \tn}-Y_{n}(t).$$

\textbf{Lemma 3.6}. {\sl Let $Y_{n}(t)$ be defined as above. \enke{Under}
the conditions of Theorem 2.2} \enke{for} any $\varepsilon>0$
 \BQN \label{eq3.6.1} \aH{\lim_{n \to
\infty}}\pk{|M_{\tn}(Y)-b_{\tn}|>\varepsilon r_{n}^{1/2}(\tn)}&= &0
\EQN and \BQN \label{eq3.6.2} \aH{\lim_{n \to
\infty}}\pk{|M_{\tn}(-Y)-b_{\tn}|>\varepsilon r_{n}^{1/2}(\tn)}&=&0
\EQN
\enke{are valid}.\\

\begin{proof}
 Since the proofs are similar, we only
give the proof of (\ref{eq3.6.1}).
\enk{By the assumptions}
\Tan{
$$\rho_{n,\tn}(t)=\frac{r_{n}(t)-r_{n}(\tn)}{1-r_{n}(\tn)}=1-c_{n}(\tn)|t|^{\alpha}+\epsilon_{n}(t)|t|^{\alpha}$$
as $t\rightarrow 0$, where
$c_{n}(\tn)=\frac{c_{n}}{1-r_{n}(\tn)}\rightarrow 1$, as
$n\rightarrow \infty$, and
$\epsilon_{n}(t)=\frac{\varepsilon_{n}(t)}{1-r_{n}(\tn)}\rightarrow
0$ as $t\rightarrow 0$, uniformly in $n$. \enk{Furthermore, for any}
$\varepsilon>0$, there exists $\gamma>0$ such that
$\sup\{|\rho_{n,\tn}(t)|, T\geq |t|\geq\varepsilon,
n\in\mathbb{N}\}<\gamma<1$. \enk{Utilising} the stationarity of $\{Y_{n}(t),
0\leq t\leq T(n)\}$,  Theorem A (i) and the definition of
$b_{T(n)}$, we have\peng{
\begin{eqnarray*}
\pk{M_{T(n)}(Y)-b_{T(n)}>\varepsilon r^{1/2}_{n}(T(n))}
&\leq & ([T(n)]+1) \pk{\max_{0\leq t\leq 1}Y_{n}(t)>\varepsilon r^{1/2}_{n}(T(n))+b_{T(n)}}\nonumber\\
&\leq& C_{6} ([T(n)]+1)  (\varepsilon r^{1/2}_{n}(T(n))+b_{T(n)})^{\frac{2}{\alpha}-1}e^{-\frac{1}{2}(r^{1/2}_{n}(T(n))+b_{T(n)})^{2}}\nonumber\\
&\leq& C_{6}([T(n)]+1)  (\log T(n))^{\frac{2-\alpha}{2\alpha}}e^{-\frac{1}{2}(2\log T(n) +\frac{2-\alpha}{\alpha}\log\log T(n)+2(r_{n}(T(n))\log T(n))^{1/2})}\nonumber\\
&\leq& C_{6}e^{-( r_{n}(T(n))\log T(n))^{1/2}}.
\end{eqnarray*}}
\enk{Assumption}
(B1) \enk{and the fact that}  $\lim_{n\to \infty} r_{n}(T(n))\log T(n)=\infty$ 
\enk{imply}
$$\aH{\lim_{n\to \infty}} \pk{M_{\tn}(Y)-b_{\tn}>\varepsilon r_{n}^{1/2}(\tn)}=0. $$}
Next, repeating the proof of equation (3.9) in Mital and Ylvisaker (1975), we have
$$\aH{\lim_{n\to \infty}} \pk{M_{\tn}(Y)-b_{\tn}<-\varepsilon r_{n}^{1/2}(\tn)}= 0,$$
hence (\ref{eq3.6.1}) \enk{holds}, and thus the claim follows. 
\end{proof}

{\it Proof of Theorem 2.2}. Represent $X_{n}(t)$ \enk{as}
$$X_{n}(t)=(1-r_{n}(\tn ))^{1/2}Y_{n}(t)+r_{n}^{1/2}(\tn)\kal{W},$$
where \aE{$\kal{W}$} is a standard \aE{Gaussian random} variable
independent of the process $\{Y_{n}(t), t\ge 0\}$. Using Lemma 3.6
and setting $a(n):= \sqrt{\frac{1-r_{n}(\tn)}{r_{n}(\tn)}}$ we \enk{obtain}
\begin{eqnarray*}
&&\pk{ r_{n}^{-1/2}(\tn)\left(\sup_{t\in[0,\tn]}|X_{n}(t)|-(1-r_{n}(\tn))^{1/2}b_{\tn}\right)\leq x}\\
&&=\pk{\sup_{t\in[0,\tn]}|X_{n}(t)|\leq r_{n}^{1/2}(\tn)[ \aH{a(n)}b_{\tn}+x]}\\
&&=\pk{-x\leq a(n) (Y_{n}(t)+b_{\tn})+\aE{\kal{W}}, a(n)(Y_{n}(t)-b_{\tn})+\aE{\kal{W}}\leq x, t\in[0,\tn]}\\
&&=\pk{a(n)(-Y_{n}(t)-b_{\tn})- \aE{\kal{W}}\leq x, a(n)(Y_{n}(t)-b_{\tn})+ \aE{\kal{W}}\leq x, t\in[0,\tn]}\\
&&=\pk{a(n)(M_{\tn}(-Y)-b_{\tn})-\aE{\kal{W}}\leq x, a(n)(M_{\tn}(Y)-b_{\tn})+\aE{\kal{W}}\leq x}\\
&&\rightarrow \pk{-\aE{\kal{W}}\leq x, \aE{\kal{W}}\leq x}, \quad n \to \infty,
\end{eqnarray*}
and hence the claim follows. \QED

{\bf Acknowledgment.}
 We would like to thank the referees and the
editor for their comments and suggestions which greatly improved the
manuscript. Z.\ Tan  has been supported by the National Science Foundation of China 11071182,
E.\ Hashorva has been supported by the Swiss National Science Foundation Grant \enke{200021-1401633/1},
Z. Peng has been supported by the National Natural Science Foundation of China 11171275.


\begin{thebibliography}{100} \small

\bibitem{}
Adler, R.J., 1990. {\it An Introduction to Continuity, Extrema, and Related Topics for General
Gaussian Processes}, Inst. Math. Statist. Lecture Notes Monogr. Ser. 12, Inst. Math.
Statist., Hayward, CA.


\bibitem{} Albin, J.M.P., Choi, H., A new proof of an old result by Pickands. Elect. Comm. in Probab., (2010), 15: 339-345.


\bibitem{} Berman, M.S., Sojourns and Extremes of Stochastic Processes, Wadsworth \& Brooks/ Cole, Boston, 1992.



\bibitem{}
\aE{D\c{e}bicki, K., Kisowski, P., A note on upper estimates for Pickands constants. Stat. Prob. Letters, 2009, 78: 2046-2051.}

\bibitem{} \aE{Dieker, A.B., Extremes of Gaussian processes over an infinite horizon. Stochastic Process. Appl., 2005, 115: 207-248.}


\bibitem{} Durrett, R., Probability theory and examples, Duxbury press, Boston, 2004.



\bibitem{} Ho, H.C., McCormick, W.P., Asymptotic distribution of sum and maximum for Gaussian processes, J. Appl. Probab., 1999, 36, 1031-1044.

\bibitem{} H\"{u}sler, J., Piterbarg, V.I., Seleznjev, O.V., On convergence of the uniform norms for Gaussian
processes and linear approximation problems. Ann. Appl. Probab. 2003, 13: 1615-1653.

\bibitem{} H\"{u}sler, J., Extremes of Gaussian processes, on results of Piterbarg and Seleznjev. Statist. Probab. Lett. 1999, 44: 251-258.


\bibitem{} Leadbetter, M.R., Lindgren, G., Rootz\'{e}n, H., Extremes and Related Properties of Random Sequences and Processes. Series in Statistics, Springer, New York, 1983.

\bibitem{} Mittal, Y., Ylvisaker, D.,  Limit distribution for the maximum of stationary Gaussian processes.
Stochastic. Process. Appl., 1975, 3: 1-18.


\bibitem{PIC}
\enk{Pickands, J. III., Asymptotic properties of the maximum in a stationary Gaussian process. Transactions
of the American Mathematical Society, 1969, 145: 75-86.}


\bibitem{} \enk{
Piterbarg, V.,  On the paper by J. Pickands "Upcrosssing probabilities for stationary Gaussian processes". Vestnik Moscow. Univ. Ser. I Mat. Mekh. 27, 25-30. English transl. in Moscow Univ. Math. Bull., 1972, 27.}

\bibitem{} Piterbarg, V.I., Asymptotic Methods in the Theory of Gaussian Processes and Fields, AMS, Providence, 1996.


\bibitem{} \aE{Shao, Q.,  Bounds and estimators of a basic constant in extreme value theory of
Gaussian processes. Statistica Sinica, 1996, 6: 245-257.}


\bibitem{} Seleznjev, O.V.,  Limit theorems for maxima and crossings of a sequence
of Gaussian processes and approximation of random processes. J. Appl. Probab., 1991, 28: 17-32.

\bibitem{} Seleznjev, O.V.,   Large deviations in the piecewise linear approximation of Gaussian processes with stationary increments.
Adv. Appl. Prob., 1996, 28: 481-499.


\bibitem{} Seleznjev, O.V.,  Asymptotic behavior of mean uniform norms for sequences of Gaussian processes and fields
Extremes., 2006, 8: 161-169.



\bibitem{} Stamatovic, B., Stamatovic, S., Cox limit theorem for large excursions of a norm of Gaussian vector process. Statist. Probab. Lett., 2010, 80: 1479-1485.




\end{thebibliography}
\end{document}